\documentclass{article}
\usepackage{amsfonts}
\newcommand{\lbl}[1]{{\tt \small [#1] \,\, }\label{#1}}
\renewcommand{\lbl}[1]{\label{#1}}
\newcommand{\rf}[1]{(\ref{#1})}

\normalbaselines
       \newtheorem{Theorem}{Theorem}[section]
       \newtheorem{Proposition}{Proposition}[section]
       \newtheorem{Corollary}{Corollary}[section]
       
       \newtheorem{Lemma}{Lemma}[section]
       \newtheorem{Remark}{Remark}[section]

\newcommand{\QED}{{\fbox{}\newline }}
\newcommand{\LL}{{\mathcal L}}
\newcommand{\FF}{{\mathcal  F}}
\newcommand{\GG}{{\mathcal  G}}
\newcommand{\EE}{{\mathcal  E}}
\newcommand{\ZZ}{\mathbb{Z}}
\newcommand{\RR}{\mathbb{R}}

\newcommand{\be}{\begin{equation}}
\newcommand{\ee}{\end{equation}}
\newcommand{\conj}{\widetilde{=}}

\newcommand{\qed}{\QED}

\newenvironment{proof}{\noindent {\bf Proof.\/}}{\qed\vskip 0.1in}
\newenvironment{proofof}[1]{\noindent {\bf Proof of #1.\/}}{\qed\vskip 0.1in}
\title{Stationary random fields with linear regressions\thanks{
Running head: Random fields with linear regressions }}

\author{Wlodzimierz  Bryc\thanks{\noindent Supported in part by Taft Memorial Grant.
\newline
{\bf AMS (1991) Subject Classification:}  Primary: 60E99; Secondary: 60G10
\newline
{\tt http://ucms02.csm.uc.edu/preprint/}}
\\ 
 Department of Mathematics \\
University of Cincinnati\\
PO Box 210025\\
Cincinnati, OH 45221--0025\\
Wlodzimierz.Bryc@UC.edu 
}
\date{
June 11, 1998
 \\ Revised: 
July 27, 2000
}
\begin{document}
\maketitle 
\centerline{
This paper is dedicated to the memory of my friend and collaborator
W{\l}odzimierz Henryk Smole\'nski (1952-1998).}

\begin{abstract}
We analyze and identify stationary fields
with linear regressions and quadratic conditional variances.
  We give sufficient conditions to determine one dimensional 
distributions uniquely as normal, and as certain compactly-supported distributions.
Our technique relies on
orthogonal polynomials, which under our assumptions 
turn out to be a version of the so called continuous
 $q$-Hermite polynomials.
 
\end{abstract}
\section{Introduction}
The literature suggests that the distribution of a random sequence is often determined
uniquely by the first two conditional
moments. A classical result of P. Levy has been augmented by a number of results that 
replace assumptions on trajectories by more conditionings. 
Building on the work of Pluci\'nska
\cite{Plucinska},
Weso{\l}owski \cite{Wesolowski93} determined all stochastic processes with quadratic
conditional variances under the assumptions that lead to processes with independent increments.
Szab{\l}owski \cite{Szablowski89}, \cite{Szablowski90} shows that 
 first two conditional moments determine  the marginal 
distribution  in smooth cases.
In \cite{Bryc-Plucinska} we 
show that the  first two conditional moments characterize gaussian sequences.
For more references and additional motivation the reader is referred to \cite[Chapter 8]{Bryc95}.

In this paper we analyze   stationary sequences under the assumptions which are
 invariant with respect
to the ordering of the index set. 
Let $(X_k)_{k\in \ZZ}$ be  a square-integrable random field.
We are interested in the fields with the first two
conditional moments   given by
\begin{equation} \lbl{(A)}
E(X_k| \dots,
X_{k-2},X_{k-1},X_{k+1},X_{k+2},\dots)=L(X_{k-1}, X_{k+1})  
\end{equation}

\begin{equation} \lbl{(B)}
E(X_k^2| \dots,
X_{k-2},X_{k-1},X_{k+1},X_{k+2},\dots)=Q(X_{k-1}, X_{k+1})
\end{equation}
for all $k\in\ZZ$ , where $L(x,y)=a(x+y)+b$ is a symmetric linear polynomial, 
and 
\be \lbl{QQQ}
Q(x,y)=A(x^2+y^2)+B xy + C +D(x+y)
\ee is a symmetric quadratic
 polynomial.

 It turns out that conditions (\ref{(A)}), (\ref{(B)}) together  with some technical
assumptions
imply that the
eigenfunctions of the shifted conditional expectation are certain
$q$-orthogonal polynomials, see \rf{Main recurrence}.
We explore this connection to identify the  one-dimensional distributions
of $(X_k)$ when the coefficients of the quadratic form \rf{QQQ}
 depend on  one  real parameter $R$  in addition to the correlation coefficient
 $\rho$. 

We show that
the following 
change occurs  as the values of the parameter change.

\begin{itemize}
\item  For a range of values between  the two ``critical values" of the 
parameter random field $(X_k)$ is bounded, its one-dimensional 
distribution is uniquely determined, and has a density.
\item At the upper ``critical value" of the parameter $R$
 the distribution $\LL(X_0)$ is normal.
\item At the lower ``critical value" of the parameter $R$
 the distribution $\LL(X_0)$ is two-valued.
\end{itemize}
%

The paper is organized as follows.
Section \ref{Assumptions} lists the assumptions used throughout the paper.
Section \ref{Theorems} contains the statements of the main results. 
The next four sections contain lemmas and proofs.
In
section \ref{Markov}  we construct stationary Markov processes which
satisfy (\ref{(A)}) and (\ref{(B)}). Section \ref{Concluding} collects
concluding remarks.

\section{Assumptions}\lbl{Assumptions}
The following general assumptions will be in place throughout the paper. 
\begin{itemize}


\item 
We assume that $(X_k)$ is $L_2$-stationary and that all one dimensional distributions
are equal, 
\be \lbl{D}
X_k \conj  X_0 \hbox{ for all } k\in \ZZ
\ee
\item  Since  (\ref{(A)}) and (\ref{(B)})  are
not affected by non-degenerate linear transformations of $X_k$, 
without loss of generality we assume
$E(X_k)=0$ and $E(X^2_k)=1$. All specific assumptions about the coefficients of  $Q(x,y)$ refer to this case.

\item We put  the following restrictions on
 the values of the first two correlation coefficients
$\rho=corr(X_0,X_1)$ and $r_2=corr(X_0,X_2)$. 

\be \lbl{r2}
\rho \ne 0
\ee  

\be \lbl{r3}
r_2+1-2|\rho|>0
\ee
\end{itemize}
\begin{Remark} Note that (\ref{r2})  excludes i.~i.~d. sequences, and  (\ref{r3}) implies
$|\rho|<1$ thus excluding sequences formed by repeating a single random variable. Clearly,
both i.~i.~d. sequences and constant sequences with second moments   satisfy (\ref{(A)}) and
 (\ref{(B)}) but can have arbitrary distribution. It is somewhat surprising that for other
values of $\rho$ the distribution of $X_0$ is determined uniquely.   
\end{Remark}

\section{Theorems}\lbl{Theorems}

Under mild assumptions condition 
\rf{(A)} determines the covariances of $(X_k)$ up to a multiplicative factor.  

\begin{Theorem} \lbl{Theorem A}
Let $(X_k)$ be  an $L_2$-stationary
random field such that
(\ref{(A)}) holds true for all $k$, and  conditions \rf{r2} and \rf{r3} are satisfied. 
 Let $r_k=corr(X_0,X_k)$ denote
the correlation coefficients. 

Then
\begin{itemize}
\item[(i)] $r_k=\rho^{|k|}$ for all $k\in \ZZ$.
\item[(ii)] One-sided regressions are also linear,
\begin{equation} \label{A0}
E(X_k|\dots, X_{-1}, X_0)=r_k X_0, \,k=1,2\dots
\end{equation}

 \end{itemize}
\end{Theorem}
By symmetry,  Theorem \ref{Theorem A} implies  the following.
\begin{Corollary} Under the assumptions of Theorem \ref{Theorem A},
$E(X_k|X_0)=\rho^{|k|} X_0$ for all $k\in \ZZ$.
\end{Corollary}

We will be interested in how the properties of one dimensional distribution
 $\LL(X_0)$ change with the
values of the coefficients of the quadratic form $Q$.  

\begin{Proposition}\lbl{Prop}
Let $(X_k)$ be 
an $L_2$-stationary
standardized random field with   conditional moments given by
(\ref{(A)}) and (\ref{(B)}), and  such that conditions (\ref{r2}-\ref{r3}) are
satisfied.
Suppose that \rf{QQQ} holds with $A<1/(1+\rho^2)$. 
If $A(\rho^2+\frac{1}{\rho^2})+B<1$ then
$X_0$ is bounded.  If $A(\rho^2+\frac{1}{\rho^2})+B=1$ and $D\ne 0$  then
either $X_0$ is bounded from below ($D>0$) or from above ($D<0$).
\end{Proposition}
Proposition \ref{Prop} suggests that
the case when \be \lbl{Q}
A(\rho^2+\frac{1}{\rho^2})+B=1 \ee and $D= 0$ is of  interest. 
This case
includes the Gaussian fields, where $A=\rho^2/(1+\rho^2)^2$ and 
$B=2\rho^2/(1+\rho^2)^2$; in  fact,  we show in \cite{Bryc99} that for a Markov chain
which satisfies \rf{(A)} and \rf{(B)}, 
the coefficients  of   \rf{QQQ} either satisfy constraint (\ref{Q}), 
or $2 A + B\rho^2=1$, and that $D=0$.
Since by taking the expected value of (\ref{(B)})
we get a trivial relation
\be \lbl{TTT}
C=1-2A-B\rho^2
\ee 
(here we used $r_2=\rho^2$ by Theorem \ref{Theorem A})
this leaves just one free parameter besides $\rho$ 
on the right hand side of \rf{QQQ}. 
As this free parameter we will use a scaled version of $B$ defined by  
\be\lbl{R}
R=B(\rho+\frac{1}{\rho})^2\ee

 \begin{Theorem}\lbl{Theorem B} Let $(X_k)$ be 
an $L_2$-stationary
standardized random field with  conditional moments given by
(\ref{(A)}) and (\ref{(B)}), and  such that conditions (\ref{D}-\ref{r3}) are
satisfied. Furthermore, assume  that the coefficients on the
right hand side
of \rf{QQQ} are such that $D=0$ and \rf{Q} holds true. Then the following three statements
hold. 
\begin{itemize}

 \item[(i)]  For $0\leq R<2$ the one-dimensional  distribution
$\LL(X_0)$  has the uniquely determined symmetric distribution   supported on a
finite interval. When $R=0$,  $X_0=\pm c$ has  two values.

\item[(ii)] For $R=2$ the  one-dimensional  distribution $\LL(X_0)$ is normal.

\item[(iii)] For $R>2$  the  one-dimensional
distributions $\LL(X_0)$ are not determined by moments.

\end{itemize} 

\end{Theorem}

The infinite number of variables in the
conditional expectations in \rf{(A)} and \rf{(B)}
 cannot be easily
reduced. 

\begin{Theorem} \lbl{Theorem C} For every value of $-1<\rho<1$ there is a
non-degenerate  $L_2$-stationary random field which  satisfies condition \rf{D},
has linear regressions
\be \lbl{(A'')}
E(X_k| X_{k-1},X_{k+1})=\alpha(X_{k-1}+ X_{k+1}),  
\ee
has quadratic second moments
\be \lbl{(B'')}
E(X_k^2|X_{k-1},X_{k+1})=\alpha^2(X_{k-1}+X_{k+1})^2+C,
\ee
and \begin{itemize} \item[(i)] the
conclusion of  Theorem \ref{Theorem A}(i) fails;
\item[(ii)] the conclusion of Theorem \ref{Theorem B} fails.
\end{itemize} 

\end{Theorem}

In Section \ref{Markov},
 for every value  of parameters $0< R\leq 2$ and $0<|\rho|<1$
 we construct strictly stationary
Markov chains which satisfy the assumptions of 
Theorem \ref{Theorem B}.

\section{Proof of Theorem \protect{\ref{Theorem A}}}

We use the notation $E(\cdot|\dots,X_n)$ to denote the conditional
expectation with respect to the sigma field generated by $\{X_k:k\leq n\}$.
Occasionally, we will also write $E^{\FF}(X):=E(X| \FF)$.  

Let $r_k=E(X_0X_k)$ be the correlation coefficients, $k=0,1,\dots$, and recall that
$\rho=r_1 \ne 0$ by assumption. 

Theorem \ref{Theorem A} follows from the following two observations.

\begin{Lemma} \lbl{oldie1}
Under the assumptions of Theorem \ref{Theorem A} 
\begin{itemize}
\item[(i)] $r_k\to 0$ as $k\to \infty$.
\item[(ii)] $r_k^2 < 1$ for all $k>2$. 
\end{itemize}
\end{Lemma}

\begin{proof} Multiplying  \rf{(A)} by $X_0$ and taking the expected value of both sides we
get $r_k=a(r_{k-1}+r_{k+1})$. Since $\rho\ne 0$,  \rf{r3} implies that $1+r_2>0$. Therefore
$a=\rho/(1+r_2)$ and  the correlation coefficients $r_k$ satisfy the recurrence $$
(1+r_2)r_k=\rho(r_{k-1}+r_{k+1}), k=1,2,\dots $$ 
For fixed $\rho,r_2$, this is a linear recurrence for $(r_k)_{k>2}$. Inequality \rf{r3}  implies
that the characteristic equation of the recurrence has two distinct  real roots and their
product is $1$. There is therefore only one root, $c$,  in the interval $(-1,1)$  and 
$r_k=b c^k$. This shows that $r_k=r_2 c^{k-2}, k=2,3,\dots$. In particular, $r_k\to 0$,
and $r_k^2<1$ for $k>2$.

\end{proof}

\begin{Lemma}\lbl{oldie2}
If $(X_k)$ satisfies the assumptions of Theorem \ref{Theorem A}, then
\be\lbl{LR0} 
E(X_{n+1}|\dots, X_n)=\rho X_n
\ee
\end{Lemma}
\begin{proof} 
We first show 
that for all $n\in\ZZ, k\geq 2 $ there are coefficients
$a(k)\ne 0$ and $b(k)\in\RR$ such that
\be \label{EE1}
E(X_{n+1}| \dots, X_{n-1}, X_n,X_{n+k},X_{n+k+1},\dots)=
a(k)X_n+b(k)X_{n+k} \hbox{}
\ee
and
\be\lbl{EE1r}
E(X_{n+k-1}| \dots, X_{n-1}, X_n,X_{n+k},X_{n+k+1},\dots)=
b(k)X_n+a(k)X_{n+k}.\ee


We prove this by induction with respect to $k\geq 2$.

For $k=2$, this follows from (\ref{(A)}) with $a(2)=\frac{\rho}{1+r_2}\ne 0$, see \rf{r2}.  

For a given  $k\geq 2$, suppose that $a(k)\ne 0$ and both
(\ref{EE1}) and \rf{EE1r} hold true for all 
$n\in\ZZ$. We will prove that the same statement holds true for $k+1$.

Conditioning on additional variable
$X_{n+k}$ and using the induction assumption we get $$E(X_{n+1}|\dots, X_{n-1},
X_n,X_{n+k+1},X_{n+k+2},\dots)=$$ 
$$E^{\dots,
X_{n-1}, X_n,X_{n+k+1},\dots}\left(E^{\dots
X_{n-1}, X_n,X_{n+k},X_{n+k+1},\dots}(X_{n+1})\right)$$
$$a(k)X_n+b(k)E(X_{n+k}|\dots, X_{n-1}, X_n,X_{n+k+1},\dots)$$
Now adding $X_{n+1}$ 
to the condition in $E(X_{n+k}|\dots, X_{n-1}, X_n,X_{n+k+1},\dots)$ 
and using the induction assumption
we get 
$$E(X_{n+k}|\dots, X_{n-1},
X_{n},X_{n+k+1},X_{n+k+2},\dots)=$$ 
$$E^{\dots, X_{n-1},
X_n,X_{n+k+1},X_{n+k+2},\dots}\left(E^{\dots, X_{n},
X_{n+1},X_{n+k+1},X_{n+k+2},
\dots}(X_{n+k})\right)=$$       
$$b(k)E(X_{n+1}|\dots, X_{n-1},
X_n,X_{n+k+1},X_{n+k+2},\dots)+a(k)X_{n+k+1}$$

Combining these two expressions we get the linear equation 
$$
M=a(k)X_{n}+b^2(k)M+a(k)b(k)X_{n+k+1}
$$
for unknown random variable $M=E(X_{n+1}| \dots, 
X_{n-1}, X_n,X_{n+k+1},X_{n+k+2},\dots)$. 
Notice that since $a(k)\ne 0$, if $b^2(k)=1$ then $X_n,X_{n+k+1}$ are linearly related. Therefore,
$r^2_{k+1}=1$, contradicting Lemma \ref{oldie1}(ii). (Case $k=2$ requires separate
verification:  $b^2(2)=\frac{\rho^2}{(1+r_2)^2}< 1$ by \rf{r3}).
Thus $b^2(k)\ne 1$ and $M$ is determined uniquely as the linear function of 
$X_n,X_{n+k+1}$. 
This proves \rf{EE1}. By symmetry (or by similar reasoning) \rf{EE1r} holds true.
The equation gives $a(k+1)=\frac{a(k)}{1-b^2(k)}$ which shows that $a(k+1)\ne 0$.
This completes the induction proof of \rf{EE1}.

Since $|r_k|<1$ for $k>2$, the
coefficients $a(k), b(k)$ in a linear regression are  determined  uniquely.
A calculation gives 
 $b(k)=\frac{r_{k-1}-r_{1}r_k}{1-r_k^2}$
Using Lemma \ref{oldie1}(i) we have $b(k)\to 0$ as $k\to\infty$. 
Passing to the limit as $k\to\infty$ in 
(\ref{EE1}) we get \rf{LR0}.
\end{proof}

\begin{proofof}{Theorem \protect{\ref{Theorem A}}}
Applying $k$-times formula \rf{LR0},  we get $E(X_k|\dots,X_0)=\rho^kX_0$.
This implies $E(X_0X_k)=\rho^k$.
\end{proofof}
\section{Proof of Proposition \protect{\ref{Prop}}}

A simple calculation using \rf{Q} and \rf{R} gives the following. 
\begin{Lemma} \lbl{L5}
The following two conditions are equivalent:
\begin{itemize}
\item[(i)]$A< 1/(1+\rho^2)$;
\item[(ii)] $R>1-1/\rho^4$.
\end{itemize}
\end{Lemma}

The second part of the next lemma will serve as the first step
in  the induction proof of
Lemma \ref{eigenfunction}.
\begin{Lemma}\lbl{E0}
If $(X_k)$ satisfies the assumptions of Theorem \ref{Theorem A} and \rf{(B)} holds true, then
\be \lbl{OldStuff}
(1-A(1+\rho^2))E(X_1^2|\dots, X_0)=(A(1-\rho^2)+B\rho^2)X_0^2 +C+D(1+\rho^2) X_0.
\ee
If in addition \rf{Q} holds true and $R\geq 0$ then
\be
\lbl{V0}
E(X_1^2|X_0)=\rho^2 X_0^2+1-\rho^2+\gamma X_0,
\ee
where $\gamma=\frac{D}{1-(1+\rho^2)A}$.
\end{Lemma}
\begin{proof}
Since $E(X_1 X_2|{\dots, X_0})=E^{\dots, X_0}\left(X_1E^{\dots, X_1}(X_2)\right)$, from 
Lemma \ref{oldie2}
 we get 
\be\lbl{miss1}
E(X_1 X_2|{\dots, X_0})=\rho E(X_1^2|\dots, X_0).
\ee

We now give another expression for the left hand side of \rf{miss1}.
Theorem \ref{Theorem A} implies that $L(x,y)=\frac{\rho}{1+\rho^2}(x+y)$.  Since
 $E(X_1 X_2|{\dots, X_0})=E(X_2E(X_1|\dots,X_0,X_2,\dots)|{\dots, X_0})$ from (\ref{(A)})
we get $E(X_1 X_2|{\dots, X_0})= \frac{\rho}{1+\rho^2}E(X_2 (X_2+X_0)|{\dots, X_0})$. 
By Lemma \ref{oldie2}
this implies $E(X_1
X_2|{\dots, X_0})=\frac{\rho^3}{1+\rho^2}X_0^2+\frac{\rho}{1+\rho^2}E(X_2^2|{\dots, X_0})$.
Since $\rho\ne 0$, combining the latter with \rf{miss1} we have
\be \lbl{OldStuff2}
E(X_2^2|{\dots, X_0})=(1+\rho^2)E(X_1^2|{\dots, X_0})-\rho^2 X_0^2.
\ee
We now substitute expression \rf{OldStuff2} into \rf{(B)} as follows. 
Taking the conditional expectation $E(\cdot|\dots,X_0)$ of both sides of
(\ref{(B)}) with $k=1$    we get
$$
E(X_1^2|{\dots, X_0})=A X_0^2+A E(X_2^2|{\dots, X_0})+B X_0^2 \rho^2+C +D(1+\rho^2) X_0.
$$
Replacing $E(X_2^2|{\dots, X_0})$ in the above expression by the right hand side of \rf{OldStuff2}
 we get \rf{OldStuff}. 

If \rf{Q} holds true then identity \rf{V0} follows by a simple calculation, since
by Lemma \ref{L5} and the assumption $R\geq 0$ we have $A< 1/(1+\rho^2)$. 
\end{proof}

By symmetry, we can switch the roles of $X_0,X_1$ in  (\ref{V0}).
Therefore the assumptions of \cite[Corollary 3]{Bryc85} are satisfied. This
implies integrability, see also \cite[Theorem 2]{Wesolowski93}
for
another proof.

\begin{Corollary}\lbl{Integrability}
Under the assumptions of Theorem \ref{Theorem A} $E (|X_k|^p)<\infty$ for all $p>1$.
\end{Corollary}

\subsection{Proof of Proposition \protect{\ref{Prop}}}
By  Theorem \ref{Theorem A}(ii)  and \rf{OldStuff}
$$\rho^2X_0^2=(E^{X_0}(X_1))^2\leq E^{X_0}(X_1^2)=\frac{A(1-\rho^2)+B\rho^2}{1-(1+\rho^2)A}X_0^2+
\frac{D}{1-(1+\rho^2)A}X_0+ \frac{C}{1-(1+\rho^2)A}.$$

To prove the first part, notice that if $X_0$ is unbounded then $\rho^2\leq \frac{A(1-\rho^2)+B\rho^2}{1-(1+\rho^2)A}$. 
The latter implies that $A(\rho^2+1/\rho^2)+B \geq 1$.
Indeed, write $A(\rho^2+1/\rho^2)+B=1+\Delta$. Then 
$\frac{A(1-\rho^2)+B\rho^2}{1-(1+\rho^2)A}=\rho^2+\frac{\rho^2\Delta}{1-(1+\rho^2)A}$, so $\Delta\geq 0$.

Furthermore, if \rf{Q} holds true then an elementary calculation gives
$\frac{A(1-\rho^2)+B\rho^2}{1-(1+\rho^2)A}=\rho^2$. Thus condition $D>0$ implies that
 $X_0$ must be bounded from below.
Similarly, if $D<0$ then $X_0$ must be bounded from above.

\section{Proof of Theorem \protect{\ref{Theorem B}}} 


\subsection{Orthogonal polynomials}

Define polynomials $Q_n$  by the recurrence
\be \lbl{Main recurrence}
xQ_n(x)=Q_{n+1}(x)+(1+q+\dots+q^{n-1})Q_{n-1}(x)
\ee
with $Q_{-1}(x)=0, Q_0(x)=1$. (Then $Q_1(x)=x$ and $Q_2(x)=x^2-1$).
Clearly, $Q_n$ are the Hermite polynomials when $q=1$. 

Polynomials $Q_n$ after a change of variable transform into the continuous
$q$-Hermite polynomials, see \cite{Ismail} and \cite{K-S}. Thus we can deduce
information about the measure that makes $Q_n$ orthogonal;
 this measure is
unique (with explicit representation for the density)  and has bounded support when $-1<q<1$,  normal when $q=1$, and is non-unique 
when $q>1$.  Later on we shall see that 
$Q_n$ are orthogonal with respect to
the distributions of $\LL(X_0)$.

The next lemma collects these results together. A sketch of the proof is
enclosed for completeness.
\begin{Lemma}\lbl{Uniqueness}
Let $Q_n$ be defined by (\ref{Main recurrence}).
\begin{itemize}
\item[(i)] If $q<-1$ then there is no probability measure with respect to which
 $Q_n$ are orthogonal.
\item[(ii)] If $-1\leq q<1$ then $Q_n$ are orthogonal with respect to a unique
 probability measure. This measure is symmetric on $\RR$ and has bounded support.  
\item[(iii)] If $q=1$ then $Q_n$ are orthogonal with respect to a unique probability measure
which is normal. 

\item[(iv)] If $q>1$ then $Q_n$ are orthogonal with respect to an infinite number
of probability
measures.
\end{itemize}
\end{Lemma}
\begin{proof}
(i) Follows from the fact that the
coefficients in  (\ref{Main recurrence}) must be non-negative.

(ii) This is an immediate consequence of Carleman's criterion,
see \cite[page 59]{Shohat-Tamarkin}. 

 It is easy to see that $q=-1$
corresponds to the (unique) symmetric
two-point distribution $X_0=\pm 1$. 
This is a degenerate case, often excluded from the general theory of orthogonal
polynomials because $Q_n(X_0)=0$ for all $n>1$.

(iii) This is the consequence of the classical result that normal distribution is
determined uniquely by its moments.

(iv) This follows from Berezanskii's result as reported in addendum 5 (page 26)
of \cite{Akhiezer}. Notice that the latter deals with normalized rather than monic
polynomials,  see (\ref{N-rec}) below. Since  $E(Q_n^2(X_0))=\prod_{k=1}^n
\frac{q^k-1}{q-1}$,  Berezanskii's theorem is used with
$b_n=\sqrt{\frac{q^{n+1}-1}{q-1}}$ and the condition $b_{n+1}b_{n-1}\leq b_n^2$
holds true because  $x\mapsto \log (q^x-1)$ is concave when $q>1$.

Several explicit weight functions for $q$-Hermite polynomials with $q>1$ 
 are given in \cite{Askey}.
\end{proof}

The following lemma  shows
that we can find orthogonal
polynomials by finding the eigenfunctions of conditional expectations, and
is based on well known motives.
\begin{Lemma}\lbl{Orthogonality}
Let $f,g$ be two functions such that $f(X_1)$ and $g(X_1)$ are square-integrable.
If $\LL(X_0)=\LL(X_1)$ and $\alpha \ne \beta$ are real numbers
such that
$E(f(X_1)|X_0)=\alpha f(X_0)$ and $E(g(X_0)|X_1)=\beta g(X_1)$
then $E(f(X_0)g(X_0))=0$.
\end{Lemma}

\begin{proof}
This follows from $E(f(X_1)g(X_0))=\alpha E(f(X_0)g(X_0))=\beta E(f(X_1)g(X_1))$.
\end{proof}
\subsection{Proof of Theorem {\protect\ref{Theorem B}}}

Let 
\be\lbl{q}
q=\frac{\rho^4+R-1}{1+\rho^4(R-1)}.
\ee

It is easy to see that  the range $0\leq R\leq 2$ 
corresponds to $-1\leq q \leq 1$ and that  $R=2$ when $q=1$, and 
 $R=0$ when $q=-1$.

\begin{Lemma}\lbl{eigenfunction}
If $q$ is defined by (\ref{q}) then polynomials $Q_n(x)$ defined by
(\ref{Main recurrence}) satisfy
\be \lbl{**}
E^{X_0}Q_n(X_1)=\rho^n Q_n(X_0).
\ee
\end{Lemma}

\begin{proof} The proof is by mathematical induction with respect to $n$.
The result is trivially true for $n=-1$ (with $Q_{-1}=0$) and for $n=0$.
By Theorem \ref{Theorem A},  (\ref{**}) holds true for $n=1$.
Since we assume that $D=0$, formula \rf{V0} implies that   (\ref{**}) holds true for $n=2$. 

Suppose (\ref{**}) holds true for $n, n-1, n-2$, and $n-3$. 
We will show that
 (\ref{**}) holds for $n+1$.

The proof repeatedly uses conditional expectations
 $\EE(\cdot)=E^{\dots,X_{-1},X_0}(\cdot)$.

Consider $\EE(X_1Q_{n}(X_2))$. We will first obtain two
equations (\ref{A'}), (\ref{B'}) which allow us to determine conditional moments 
$\EE(X_1Q_{n}(X_2))$ and
$\EE(X_1Q_{n}(X_1))$.

The first equation is obtained as follows. By induction assumption,
$\EE(X_1Q_{n}(X_2))=\EE(X_1E^{\dots,X_{0},X_1}Q_{n}(X_2))=
\rho^n \EE(X_1Q_{n}(X_1))$. On the other hand, 
equation (\ref{(A)}) gives 
$$\EE(X_1Q_{n}(X_2))=
\frac{\rho}{1+\rho^2} 
(\EE(X_0Q_{n}(X_2))+\EE(X_2Q_{n}(X_2))).$$
Therefore
\be\lbl{A'}
\EE(X_2Q_{n}(X_2))=
(1+\rho^2)\rho^{n-1} \EE(X_1Q_{n}(X_1))- \rho^{2n} X_0Q_{n}(X_0).
\ee

To obtain the second equation, consider $\EE(X_1^2Q_{n-1}(X_2))$.
By induction assumption, $\EE(X_1^2Q_{n-1}(X_2))=
\rho^{n-1} \EE(X_1^2Q_{n-1}(X_1))$. On the other hand,
 equations (\ref{(B)}) and \rf{QQQ} with $D=0$  give 
$$\EE(X_1^2Q_{n-1}(X_2))=$$
$$A \EE(X_2^2Q_{n-1}(X_2))+
A \EE(X_0^2Q_{n-1}(X_2))+$$
$$B \EE(X_0X_2Q_{n-1}(X_2))+
C \EE(Q_{n-1}(X_2)).
$$
(We will use relations \rf{Q} and \rf{TTT} later on.)

Thus, using the induction assumption, we obtain the second equation
\begin{eqnarray}
\lbl{B'}
&\rho^{n-1} \EE(X_1^2Q_{n-1}(X_1))=&\\
&A \EE(X_2^2Q_{n-1}(X_2))+\nonumber
A \rho^{2n-2}  X_0^2Q_{n-1}(X_0)+&\\
&B X_0 \EE(X_2Q_{n-1}(X_2))+
C \rho^{2n-2} Q_{n-1}(X_0).&\nonumber
\end{eqnarray}

Polynomials $Q_n$ satisfy second order
 recurrence of the form $x Q_n(x)=Q_{n+1}(x)+\beta_{n+1}Q_{n-1}(x)$.
We use it to rewrite (\ref{B'}) as follows.
\begin{eqnarray}\lbl{B''} 
&\rho^{n-1} \EE(X_1Q_{n}(X_1))+\rho^{n-1}\beta_{n}\EE(X_1Q_{n-2}(X_1)) =&\\
&A \EE(X_2Q_{n}(X_2))+A\beta_n\EE(X_2Q_{n-2}(X_2))+\nonumber
A \rho^{2n-2}  X_0Q_{n}(X_0)+A \rho^{2n-2}\beta_n X_0Q_{n-2}(X_0)+&\\
&B X_0 \EE(Q_{n}(X_2))+B X_0 \beta_n\EE(Q_{n-2}(X_2))+
C \rho^{2n-2} Q_{n-1}(X_0).&\nonumber
\end{eqnarray}
Now use (\ref{A'}) 
 to eliminate $\EE(X_2Q_{n}(X_2))$ from this equation.
Since $C=1-2A-\rho^2 B$ this gives
\begin{eqnarray}\lbl{B'''} 
&(1-A(1+\rho^2))\EE(X_1Q_{n}(X_1)) =&\\
&\beta_{n}(A(1+\rho^2)\rho^{-2} -1)\EE(X_1Q_{n-2}(X_1))+\nonumber&\\
&\beta_n\rho^{n-3}(A (\rho^{2} -
 1)+B )X_0 Q_{n-2}(X_0)+\nonumber&\\
&\rho^{n-1}(A (1  - \rho^{2}) +B \rho^{2} )X_0Q_{n}(X_0)+&\nonumber\\
&(1-2A-\rho^2 B) \rho^{n-1} Q_{n-1}(X_0).&\nonumber
\end{eqnarray}

By (\ref{Q}), identity 
$\frac{A (\rho^{2} +1)-\rho^2 }{\rho^2(1-A(1+\rho^2))}=-\rho^2 q$ 
(recall that $q$ is given by
(\ref{q})) 
 and elementary calculation we now get

\begin{eqnarray}\lbl{B4} 
&\EE(X_1Q_{n}(X_1)) =&\\
&\beta_{n}q\rho^2\left(\rho^{n-3} X_0 Q_{n-2}(X_0)-\EE(X_1Q_{n-2}(X_1))\right)+&\nonumber\\
&\rho^{n+1}X_0Q_{n}(X_0)+
(1-\rho^2) \rho^{n-1} Q_{n-1}(X_0)&\nonumber
\end{eqnarray}
Since $x Q_n(x)=Q_{n+1}(x)+\beta_{n+1}Q_{n-1}(x)$, and 
by induction assumption 
$\EE(Q_{n-3}(X_1))=\rho^{n-3}Q_{n-3}(X_0)$, this implies

\begin{eqnarray*}
&\EE(Q_{n+1}(X_1))=\rho^{n+1}Q_{n+1}(X_0)+&\\
&
(1-\rho^2)\rho^{n-1}\left(q\beta_n-
\beta_{n+1}+1 
\right) Q_{n-1}(X_0).&
\end{eqnarray*}

Since (\ref{Main recurrence}) means that $\beta_{n+1}=q\beta_n + 1$, this
 ends the proof.  
\end{proof}

\begin{proofof}{Theorem {\protect\ref{Theorem B}}} 
All assumptions are symmetric with respect to time-reversal. Therefore, formula (\ref{**})
implies  $E(Q_n(X_0)|X_1)=\rho^n
Q_n(X_1)$. Since $0<|\rho|<1$,
 Lemma \ref{Orthogonality} (used with $f=Q_m, g=Q_n, n\ne m$) proves that $Q_n$
are  orthogonal with respect to $\LL(X_0)$. Therefore,  Lemma
\ref{Uniqueness}  identifies uniquely
 the distribution $\LL(X_0)$ for both
$R=2$ and $R<2$ cases. 

The distribution is not determined uniquely by the
moments when $q>1$ which corresponds to $R>2$. 

Finally, when the distribution is determined uniquely, the odd-order moments $EX_0^{2n+1}=0$ by
(\ref{Main recurrence}). Therefore the distribution of $X_0$ is symmetric.

%
%
%
%
%
\end{proofof}
\section{Proof of Theorem \protect{\ref{Theorem
C}}}

We begin with the following simple lemma.
\begin{Lemma}\lbl{indep}
Let $X_0,Y_0$ be integrable random variables.
Suppose that $\FF,\GG$ are $\sigma$-fields such that 
 $\sigma(X_0,\FF)$ and $\sigma(Y_0,\GG)$ are
 independent, 
and  
there is $\rho$ such that
$
E(X_0|\FF)=\rho X_1
$, 
$
E(Y_0|\GG)=\rho Y_1
$,
$
E(X_0^2|\FF)=\rho^2 X_1^2+1-\rho^2
$
, and $
E(Y_0^2|\GG)=\rho^2 Y_1^2+1-\rho^2
$.

Let $Z_k=(a X_k+b Y_k)$, where $a^2+b^2=1$ and denote by 
$\FF\vee\GG$ the $\sigma$-field generated by $\FF\cup\GG$ .
 Then
\be \lbl{Z1}
E(Z_0|\FF\vee\GG)=\rho Z_1,
\ee
\be \lbl{Z2}
E(Z_0^2|\FF\vee\GG)=\rho^2 Z_1^2+1-\rho^2.
\ee
\end{Lemma}
\begin{proof}
Clearly $E(Z_0|\FF\vee\GG)=aE(X_0|\FF\vee\GG)+bE(Y_0|\FF\vee\GG)$ proving
(\ref{Z1}).
Similarly,
$E(Z_0^2|\FF\vee\GG)=a^2E(X_0^2|\FF\vee\GG)+b^2E(Y_0^2|\FF\vee\GG)+2abE(Y_0
X_0|\FF\vee\GG)$. 
Now  $E(Y_0 X_0|\FF\vee\GG)=E(Y_0 E(X_0|\FF\vee\GG,Y_0)|\FF\vee\GG)=\rho
E(Y_0 X_1|\FF\vee\GG)$. Since $X_1$ is $\FF\vee\GG$-measurable,  we get $E(Y_0
X_0|\FF\vee\GG)=\rho^2 X_1 Y_1$, which proves (\ref{Z2}).

\end{proof}

\subsection{Proof of Theorem \protect{\ref{Theorem C}}} 
Fix $-1<\rho<1$. 
First define a periodic stationary sequence $(\xi_k)$ such that
 $\xi_{k+2}=\xi_k$ and with the correlation coefficient $\rho$.
To this end define the joint distribution of $\xi_1,\xi_2$ by

$$ \Pr(\xi_1=1,\xi_2=1)=\Pr(\xi_1=-1,\xi_2=-1)=\frac{1+\rho}{4},$$
$$ \Pr(\xi_1=-1,\xi_2=1)=\Pr(\xi_1=1,\xi_2=-1)=\frac{1-\rho}{4}.$$

Let $(\gamma_k)$ be a centered  Markov gaussian sequence with correlations
$E(\gamma_0\gamma_k)=r^k$, where $r=\frac{1-\sqrt{1-\rho^2}}{\rho}$, and
independent of $(\xi_k)$.

Let $X_k=a \xi_k+ b \gamma_k$, 
where $a^2+b^2=1$. 
 Then $E(X_0X_1)=a^2\rho +b^2 r$.
By selecting $a$ close enough to $1$, and by varying $\rho$ we can thus
have correlations $corr(X_0,X_1)$ fill out the entire interval $(-1,1)$.

Using Lemma \ref{indep} we verify that $(X_k)$ satisfies (\ref{(A'')})
and (\ref{(B'')}) with $\alpha=\rho/2, C=1-\rho^2$. Indeed, $E(\xi_1|\xi_0,\xi_2)=\rho\xi_0=\rho(\xi_0+\xi_2)/2$
and $E(\xi_1^2|\xi_0,\xi_2)=1=\rho^2(\xi_0+\xi_2)^2/4+1-\rho^2$.
Similarly, the gaussian sequence satisfies
$E(\gamma_1|\gamma_0,\gamma_2)=\frac{r}{1+r^2}(\gamma_0+\gamma_2)=\rho(\gamma_0+\gamma_2)/2$
and $E(\gamma_1^2|\gamma_0,\gamma_2)=\rho^2(\gamma_0+\gamma_2)^2/4+1-\rho^2$

When $a\ne 0$ the correlation coefficients $E(X_0X_k)$ do not
converge to $0$. Thus the conclusion of Theorem \ref{Theorem A}(ii) fails. 

Notice that $Q=L^2+const$ which corresponds to the normal case 
($R=0$); 
the distribution of $X_0$ is not compactly-supported 
  but for $a\ne 0$  it is not normal. Thus the conclusion of
Theorem \ref{Theorem B} fails. 

\section{Construction of Markovian fields} \lbl{Markov}
In this section we construct Markov fields which satisfy
(\ref{(A)}) and (\ref{(B)}), proving the following.

\begin{Proposition}\lbl{TM}
For all $0<\rho^2<1, 0\leq  R\leq 2$,
there are stationary Markov chains that satisfy the assumptions 
of Theorem \ref{Theorem B}. 
\end{Proposition}

Since $R=0$  corresponds to the elementary two-valued Markov chain 
which was explicitly analyzed in \cite{Bryc99}, and $R=2$ corresponds to 
stationary Gaussian Markov processes,
 we  only   consider the case $0<R<2$. 

Through the remainder of this section we fix $\rho, R$.
Let $Q_n$ be given by the recurrence 
\be\lbl{N-rec}
xQ_n(x)=b_{n+1}Q_{n+1}(x)+b_nQ_{n-1}(x) 
\ee
with initial polynomials $Q_0(x)=1,Q_1(x)=x$,
where $b_n=\sqrt{1+q+\dots+q^{n-1}}, b_0=0, b_1=1$, and $q$ is given by \rf{q}. 
Since  $q>-1$, $b_n>0$ for all $n\geq 1$.
These are the normalized orthogonal polynomials from
the proof of Theorem \ref{Theorem B}.
 
Let $\mu$ be a probability measure which 
makes $Q_n$ orthogonal.

\begin{Lemma}\lbl{posit}
If $-1<q\leq 1$ and $|\rho|<1$,  then
$\sum_{n=0}^\infty \rho^n Q_n(y) Q_n(x)$ converges in 
$L_2(\mu(dx)\mu(dy))$ to a non-negative function.
\end{Lemma}
\begin{proof}
The case $q=1$ is classical
and the sum of the series  is
$
\frac{1}{\sqrt{1-\rho^2}}\exp(-\rho\frac{\rho(x^2+y^2)-2 x y
}{2(1-\rho^2)})$.     

 For $-1<q<1$ 
an explicit product representation for the series can be deduced from
the facts collected in \cite{Ismail}, see also \cite[(3.9)]{Carlitz}. Namely, 
in the notation of \cite{Carlitz}, 
 $Q_n(x)=\frac{e^{-i n\theta} H_n(e^{2i n\theta})}{(1-q)^{n/2}b_{1} \dots b_n}$ where 
$x=2 \cos (\theta)/\sqrt{1-q}$. Therefore, using
\cite[(3.9)]{Carlitz} with
$x=2 \cos (\theta_x)/\sqrt{1-q}$, $y=2 \cos (\theta_y)/\sqrt{1-q}$ (see also
\cite[(2.10), (5.2), and (5.7)]{Ismail}) we have

$$\sum_{n=0}^\infty \rho^n Q_n(y) Q_n(x)=
\prod_{k=0}^\infty \frac{(1-\rho^2 q^k)}{
(1+\rho^2q^{2k}-2\rho q^k \cos(\theta_x+\theta_y))(1+\rho^2q^{2k}-
2\rho q^k \cos(\theta_x-\theta_y))}.$$
Since the last expression is a product of positive factors, this ends
the proof when $|q|<1$.

\end{proof}

Let $(X_n)$ be a Markov chain with initial distribution $\mu$ and
 the transition probability 
$$P_x(dy)=\sum_{n=0}^\infty \rho^n Q_n(y) Q_n(x)
\mu(dy).$$


\begin{Lemma}\lbl{LinReg} $(X_n)$ is stationary, and
satisfies condition \rf{(A)}.
\end{Lemma}
\begin{proof}
Applying Fubini's theorem to the function which by Lemma \ref{posit} is
non-negative  we get
$\int P_x(A)\mu(dx)=
 \int_A  \int_R \sum_{n=0}^\infty \rho^nQ_n(y) Q_n(x)\mu(dx)\mu(dy)$.
Since $\int Q_n(x)\mu(dx)=0$ for all $n>0$ and $Q_0=1$, we get $\int P_x(A)\mu(dx)=\mu(A)$.
This proves stationarity.

To prove \rf{(A)}, by the Markov property it suffices to show that
$E(X_1|X_0,X_2)=\frac{\rho}{1+\rho^2}(X_0+X_2)$.
Let $\phi(X_0,X_2)\geq 0$ be an arbitrary bounded measurable function. 
We will verify that  
\be
\lbl{LL}
E(X_1\phi(X_0,X_2))=\frac{\rho}{1+\rho^2}E((X_0+X_2)\phi(X_0,X_2) )
.\ee

Since $\phi$ is square integrable, and $Q_{n,m}(x,y)=Q_n(x)Q_m(y)$ are
orthogonal in  $L_2(\mu(dx)\mu(dy))$ we can write
$\phi(X_0,X_2)=\phi_0(X_0,X_2)+\sum_{i,j=0}^\infty
\phi_{i,j}Q_i(X_0)Q_j(X_2)$, where $\phi_0$ is orthogonal to all
polynomials in variables $X_0,X_2$. In the case $R<2$ we have $\phi_0 =0$
because $(Q_n)$ are an orthogonal basis of $L_2(\mu)$. 

Notice that the joint distribution of $X_0,X_1,X_2$ is given by
\be\lbl{joint}
\nu(dx,dy,dz)=\sum_{n,k=0}^\infty
\rho^{n+k}Q_n(x)Q_n(y)Q_k(y)Q_k(z)\mu(dx)\mu(dy)\mu(dz)
.\ee 

From (\ref{N-rec}) and orthogonality we get
$$ \int y Q_n(y)Q_{n-1}(y)\mu(dy)=b_{n}$$ 
Moreover, $ \int y Q_n(y)Q_{k}(y)\mu(dy)=0$  for $k\ne n\pm1$.

Therefore 
$$\int y \phi(x,z) \nu(dx,dy,dz)=
\sum_{n,k} 
\rho^{n+k}\phi_{n,k}\int y Q_n(y)Q_k(y)\mu(dy) 
$$
$$=\rho\sum\rho^{2n}b_{n+1}(\phi_{n,n+1}+\phi_{n+1,n}).$$
Similar calculation gives 
$$\int x \phi(x,z) \nu(dx,dy,dz)=
\rho^2\sum_{n=0}^\infty \rho^{2n}b_{n+1}\phi_{n,n+1}+
\sum_{n=0}^\infty \rho^{2n}b_{n+1}\phi_{n+1,n}
.$$
Since by symmetry a similar formula holds true for the integral of $z$ instead of $x$,
we get $$\int (x+z) \phi(x,z) \nu(dx,dy,dz)=(1+\rho^2)
\sum_{n=0}^\infty \rho^{2n}b_{n+1}(\phi_{n,n+1}+\phi_{n+1,n}).$$
Comparing the coefficients in the expansions we verify that
$$\int y \phi(x,z) \nu(dx,dy,dz)=
\frac{\rho}{1+\rho^2}\int (x+z) \phi(x,z) \nu(dx,dy,dz),
,$$
proving
(\ref{LL}). 

\end{proof}

\begin{Lemma}\lbl{QuadrVar} 
$(X_n)$ satisfies condition \rf{(B)}
 with the coefficients in \rf{QQQ} 
determined by the equations \rf{TTT}, \rf{Q}, \rf{R}, and  $D=0$.
\end{Lemma}
\begin{proof}
By the Markov property it suffices to show that
 $E(X_1^2|X_0,X_2)=Q(X_0,X_2)$.
To this end, as in the proof of Lemma \ref{LinReg},
we fix  a  bounded measurable function $\phi(X_0,X_2)=\sum_{n,k}\phi_{n,k}Q_n(X_0)Q_k(X_2)$. 
We will verify that  
\be \lbl{QQ}
E(X_1^2\phi(X_0,X_2))=E\left(\left(A(X_0^2+X_2^2+BX_0X_2+C)\right)\phi(X_0,X_2)\right)
\ee
by comparing the coefficients in the orthogonal expansions.
 
From (\ref{N-rec}) we get
$ \int y^2 Q_n^2(y) \mu(dy)=b_{n+1}^2+b_n^2$, 
and
$ \int y^2 Q_n(y)Q_{n+2}(y)\mu(dy)=b_{n+1}b_{n+2}$.
Moreover, by orthogonality
 $ \int y^2 Q_n(y)Q_{n+k}(y)\mu(dy)=0$ except when $k=0,2,-2$. 
Using these identities and the expansion
$\int y^2 \phi(x,z)\nu(dx,dy,dz)=\sum_{n,k}\rho^{n+k}\phi_{n,k}\int y^2 Q_n(y)Q_{n+k}(y)\mu(dy)$
we see that
\begin{eqnarray}
&\int y^2 \phi(x,z)\nu(dx,dy,dz)=&\\
&\sum_{n=0}^\infty \rho^{2n}\phi_{n,n}(b_{n+1}^2+b_n^2)+
\rho^2 \sum_{n=0}^\infty \rho^{2n}b_{n+1}b_{n+2}(\phi_{n,n+2}+\phi_{n+2,n})
.&\nonumber
\end{eqnarray}
We now turn to the right hand side of \rf{QQ}.
Since $$\int x^2 \phi(x,z)\nu(dx,dy,dz)=\int \phi(x,z)\sum_{n=0}^\infty  \rho^{2n}
x^2 Q_n(x)Q_n(z)\mu(dx)\mu(dz),$$
a calculation based on \rf{N-rec} yields
\begin{eqnarray}
&\int x^2 \phi(x,z)\nu(dx,dy,dz)=&\\
&\sum_{n=0}^\infty \rho^{2n}\phi_{n,n}(b_{n+1}^2+b_n^2)+
\sum_{n=0}^\infty \rho^{2n}\phi_{n,n+2}b_{n+1}b_{n+2}+
\rho^4 \sum_{n=0}^\infty \rho^{2n}\phi_{n+2,n}b_{n+1}b_{n+2}.&\nonumber
\end{eqnarray}
By symmetry,
\begin{eqnarray}
&\int (x^2+z^2) \phi(x,z)\nu(dx,dy,dz)=&\\
&2\sum_{n=0}^\infty \rho^{2n}\phi_{n,n}(b_{n+1}^2+b_n^2)+
(1+\rho^4) \sum_{n=0}^\infty \rho^{2n}(\phi_{n,n+2}+\phi_{n+2,n})b_{n+1}b_{n+2}
.&\nonumber
\end{eqnarray}
Another elementary calculation using \rf{N-rec} gives
\begin{eqnarray}
&\int xz \phi(x,z)\nu(dx,dy,dz)=&\\
&\frac{1}{\rho^2}\sum_{n=1}^\infty \rho^{2n}\phi_{n,n}b_n^2+
\rho^2\sum_{n=0}^\infty \rho^{2n}\phi_{n,n}b_{n+1}^2+&\nonumber\\
&\rho^2 \sum_{n=0}^\infty \rho^{2n}(\phi_{n,n+2}+\phi_{n+2,n})b_{n+1}b_{n+2}.&\nonumber
\end{eqnarray}

Since $\int \phi(x,z)\nu(dx,dy,dz)=\sum_{n=0}^\infty \rho^{2n}\phi_{n,n}$,
we can now verify that \rf{QQ} holds true by verifying the relation
$$
\int y^2 \phi(x,z)\nu(dx,dy,dz)=\int \left(A(x^2+z^2)
+B xz +C\right) \phi(x,z)\nu(dx,dy,dz)
.$$
Comparing the coefficients in the expansions, the latter 
reduces to the following.
\begin{itemize}
\item Coefficients at $\phi_{0,0}$  match when 
$$ (2A+B\rho^2)b_1^2+C=b_1^2
$$
(this holds true by \rf{TTT})
\item Coefficients at $\phi_{n,n}$ for $n\geq 1$ match when
$$
2A(b_{n+1}^2+b_n^2)+B\rho^2b_{n+1}^2+B\frac{b_n^2}{\rho^2}+C=b_{n+1}^2+b_n^2
$$
(this holds true after a longer calculation using $b_{n+1}^2=qb_n^2+1$, \rf{Q}, 
and \rf{TTT})
\item Coefficients at $\phi_{n,n+2}$ and at $\phi_{n+2,n}$ for $n\geq 0$ match when
$$ A(1+\rho^4)+B\rho^2=\rho^2 
$$
(this holds true by \rf{Q})
\end{itemize}
 This implies that \rf{QQ} holds true.
\end{proof}

\section{Concluding Remarks}\lbl{Concluding}

\begin{enumerate}
\item To simplify the notation, 
we restricts ourselves to one-dimensional 
distributions of the random fields on $\ZZ$. In principle
our technique is directly applicable
to all finite-dimensional distributions, 
and to multivariate-valued random fields on $\ZZ$, see \cite{Bryc99}.

\item Symmetric
distributions in Theorem \ref{Theorem B}  differ from
the
distributions corresponding to quadratic conditional variances
 in \cite{Wesolowski93}.
\item The density in Theorem \ref{Theorem B} (i) has
 explicit
product
representation which can be recovered by a change of variable 
from the density identified as the 
``weight function" for the continuous $q$-Hermite
polynomials, see \cite[Section 3.29]{K-S} and \cite[(2.14)]{Allaway}.

 An interesting 
explicit case arises when 
$Q(x,y)=\frac{\rho^2}{1+\rho^2}\left(x^2+y^2 \right)+\frac{\rho^2(1-\rho^2)}{1+\rho^2}xy + \frac{(1-\rho^2)(1-\rho^3)}{1+\rho^2}$. 
In this case $q=0$ and thus
$Q_n(x/2)/2$  
are  the Chebyshev polynomials of the second kind. This implies that
 $|X_0|\leq 2$ has the density $\frac{1}{2\pi}\sqrt{4-x^2}$. 

\item It would be interesting to know whether 
there are stationary processes which satisfy conditions \rf{(A)}, \rf{(B)}
with coefficients in \rf{QQQ} corresponding to $R>2$. 
(Formula \cite[(3.13)]{Carlitz} with $q'=1/q<1$ indicates that
the conclusion of Lemma \ref{posit} fails.)

\end{enumerate}

\subsection*{Acknowledgement}
I would like to thank M. E. H. Ismail, 
W. Matysiak,   P. Szab{\l}owski, and the referee for helpful comments and information.

\subsubsection*{Note Added Late}
 Matysiak \cite{Matysiak99} pointed out to us that  assumption
\rf{r3} can be weakened to  $1+r_2-2\rho^2>0$;
the latter condition is
just the linear independence of $X_0,X_1,X_2$. 

\end{document}